\documentclass[12pt]{amsart}
\usepackage{amscd,amssymb}
\usepackage[graph,frame,poly,arc]{xy}  
\usepackage{tikz}
\usepackage[plainpages,backref,urlcolor=blue]{hyperref}

\theoremstyle{plain}
\newtheorem{thm}[subsection]{Theorem}
\newtheorem{lem}[subsection]{Lemma}

\newtheorem{cor}[subsection]{Corollary}

\theoremstyle{definition}
\newtheorem{rk}[subsection]{Remark}

\newtheorem{ex}[subsection]{Example}

\numberwithin{equation}{section}
\newcommand{\A}{{\mathcal A}}

\newcommand{\D}{{\mathcal D}}

\newcommand{\al}{{\alpha}}
\newcommand{\be}{{\beta}}

\newcommand{\C}{\mathbb{C}}
\newcommand{\PP}{\mathbb{P}}

\DeclareMathOperator{\rank}{rank}

\begin{document}
%\date{June 4, 2009}

\title [Waring rank of binary forms, harmonic cross-ratio and golden ratio]
{Waring rank of binary forms, harmonic cross-ratio and golden ratio}

\author[Alexandru Dimca]{Alexandru Dimca}
\address{Universit\'e C\^ ote d'Azur, CNRS, LJAD, France}
\email{dimca@unice.fr}

\author[Gabriel Sticlaru]{Gabriel Sticlaru}
\address{Faculty of Mathematics and Informatics,
Ovidius University
Bd. Mamaia 124, 900527 Constanta,
Romania}
\email{gabrielsticlaru@yahoo.com }

%\thanks{$^1$ This work has been partially supported by the French government, through the $\rm UCA^{\rm JEDI}$ Investments in the Future project managed by the National Research Agency (ANR) with the reference number ANR-15-IDEX-01 and by the Romanian Ministry of Research and Innovation, CNCS - UEFISCDI, grant PN-III-P4-ID-PCE-2016-0030, within PNCDI III.}

\subjclass[2010]{Primary 14J70; Secondary  14B05, 32S05, 32S22}

\keywords{Waring decomposition, Waring rank, 
line arrangement, cross-ratio, golden ratio}

\begin{abstract} 
We discuss the Waring rank of binary forms of degree 4 and 5, without multiple factors, and point out unexpected relations to the harmonic cross-ratio, j-invariants and the golden ratio. These computations of ranks for binary forms are used to show that the combinatorics of a line arrangement in the complex projective plane does not determine the Waring rank of the defining equation even in very simple situations.
\end{abstract}
 
\maketitle

%\tableofcontents

\section{Introduction} 
For the general question of symmetric tensor decomposition we refer to
\cite{Guide, Car1+, Ci, CGLM, CS, FOS, FLOS, IK,L1,LT,MO,O}, as well as to the extensive literature
quoted at the references in \cite{Guide} and \cite{FLOS}.
 Consider the graded polynomial ring $S=\C[x,y]$, let $S_d$ denote the vector space of homogeneous polynomials of degree $d$ in $S$, and let $f \in S_d$ be a binary form of degree $d$.  We consider the {\it Waring decomposition}
\begin{equation}
\label{eq1}
(\D) \ \  \  \ \  f=\ell_1^d + \cdots +\ell_r^d,
\end{equation}
where $\ell_j \in S_1$ are linear forms in $x,y$, and $r$ is minimal, in other words $r=\rank f$ is the {\it Waring rank}  of $f$. Hence, the nonzero binary form $f$ has Waring rank $r=\rank f=1$ if and only if $f$ is the power of a linear form. Note that the Waring rank of a form $f$ of degree $d$ depends only on the corresponding class $[f]$ in $\PP(S_d)$, and even on the corresponding $SL_2(\C)$-orbit of $[f]$ in $\PP(S_d)$. It is clear that two forms $f$ and $f'$ in $S_d$ such that
\begin{equation}
\label{eq10}
k=\rank f =\rank f' \in \{1,2\}
\end{equation}
give rise to the same $SL_2(\C)$-orbit in $\PP(S_d)$. The rank two binary forms are discussed in detail in \cite{BM}.

In this note we discuss the Waring ranks of binary quartics and binary quintics, assuming they have distinct factors.
For binary quartics the generic rank is 3. We describe precisely the quartic forms of rank 2 in terms of the {\it harmonic cross-ratio} of the corresponding 
roots in $\PP^1$, and explain why all the other binary quartics with distinct factors have rank 3, see Theorem  \ref{thm1}.
For binary quintics with distinct factors, those of rank 2 are closely related to the {\it golden ratio}. The generic binary quintics still have Waring rank 3, and there is an algebraic curve parametrizing the binary quintics with distinct factors and with rank 4, see Theorem  \ref{thm2}.

In the final section we use the previous results to show that the combinatorics of a line arrangement $\A: F(x,y,z)=0$ in $\PP^2$ does not determine the Waring rank of $F$ even in very simple situations, namely when $F(x,y,z)=zf(x,y)$, see Theorem \ref{thm3} and Example \ref{ex4}.

\medskip

We would like to thank Alessandro Oneto for kindly drawing our attention to several key results in \cite{Car2+}, and to Laura Brustenga i Moncus\' i
for useful informations concerning \cite{BM}. Computations with CoCoa \cite{Co} and Singular \cite{Sing} also played a key role in our results.

\section{Sylvester's Theorem} 

The Waring rank $r=\rank f$ can be described as follows. Let $Q=\C[X,Y]$,  where $X=\partial_x$ and $Y=\partial_y$. Then $Q$ is the ring of differential operators with constant coefficients and acts on $S$ in the obvious way. For a binary form $f \in S$, we consider the ideal
of differential operators in $Q$ killing $f$, namely
\begin{equation}
\label{eq2}
Ann(f)=\{q \in Q : q \cdot f=0\},
\end{equation}
also denoted by $f^{\perp}$ and called the apolar ideal of $f$. Note that $Ann(f)$ is a graded ideal, whose degree $s$ homogeneous component is given by $\ker [f]$, where
\begin{equation}
\label{eq3}
[f]: Q_s \to S_{d-s}
\end{equation}
is the morphism $g \mapsto g \cdot f$. The matrix of this linear map with respect to the obvious monomial bases in $Q_s$ and $ S_{d-s}$ is called the {\it catalecticant matrix} $C(f)_s$ of $f$ in degree $s$.

\begin{ex}
\label{ex1} 
As a example, if we take $d=4$ and write
$$f=a_0x^4+4a_1x^3y+6a_2x^2y^2+4a_3xy^3+a_4y^4,$$
then
\begin{center}
$$C(f)_2=12 \left(
  \begin{array}{ccccccc}
     a_0 & a_1& a_2  \\
     2a_1 &2a_2& 2a_3 \\
    a_2& a_3& a_4  \\
   \end{array}
\right).$$
\end{center}

\end{ex}

The following result is perhaps well known.
\begin{lem}
\label{lemS}
\begin{enumerate}
The graded ideal $Ann(f) \subset Q$ of the binary form $f$ of degree $d$ satisfies the following.

\item $Ann(f)_0 \ne 0$ if and only if $f=0$.

\item $Ann(f)_0 = 0$ and $Ann(f)_1 \ne 0$ if and only if $f=x^d$  after a linear change of coordinates.

\item $Ann(f)_1= 0$ and $Ann(f)_2= \C \ell^2$ for some $\ell \in Q_1$
if and only if $f=x^{d-1}y$  after a linear change of coordinates.

\end{enumerate}

\end{lem}

The following result goes back to Sylvester \cite{Sy}. See also \cite{CS}.

\begin{thm}
\label{thmS}
For a binary form $f$ of degree $d$,
the apolar ideal $Ann(f)$ is a complete intersection, namely there are two binary forms $g_1$ and $g_2$ in $Q$ such that $Ann(f)=(g_1,g_2)$.
The degrees $d_j$ of $g_j$ for $j=1,2$ satisfy $d_1+d_2=d+2$.
Moreover, if we assume $d_1 \leq d_2$, then the Waring rank $r=\rank f$ is determined as follows.
\begin{enumerate}

\item If the binary form $g_1$ has no multiple factors, then $r=d_1$.

\item Otherwise, $r=d_2$.

\end{enumerate}

\end{thm}
According to Lemma \ref{lemS}, the interesting case is 
$$2 \leq d_1 \leq d_2.$$
In this case we have the following, see also the Introduction in \cite{BT}.
\begin{thm}
\label{thmAH}
If the binary form $f$ of degree $d$ satisfies $\rank f\geq 2$, then
$$\rank f \leq d$$
and the equality holds if and only if $f=x^{d-1}y$  after a linear change of coordinates.
Moreover, for a generic binary form $f$ of degree $d$ one has
$$\rank f = \left \lfloor \frac{ d+2}{2} \right \rfloor=  \left \lceil \frac{ d+1}{2} \right \rceil.$$
\end{thm}
\proof
The first claim follows from Lemma \ref{lemS}, (3) and
Sylvester's Theorem \ref{thmS}. The second claim is a special case of 
Alexander-Hirschowitz results in \cite{AH}.
\endproof

\begin{ex}
\label{exd3} 
When $d=3$, a binary form $f$ has $\rank f=2$ if and only if $f$ has no multiple factor, and then $f$ is projectively equivalent to the binary form
$x^3+y^3$.
\end{ex}

\section{Binary quartics and the harmonic cross-ratio}
In this section we investigate the Waring rank of binary forms of degree $4$ having no multiple factor. 
If we  write
$$f=a_0x^4+4a_1x^3y+6a_2x^2y^2+4a_3xy^3+a_4y^4,$$
then
the determinant 
\begin{center}
$$ T(f)=\det \left(
  \begin{array}{ccccccc}
     a_0 & a_1& a_2  \\
     a_1 &a_2& a_3 \\
    a_2& a_3& a_4  \\
   \end{array}
\right)=a_0a_2a_4+2a_1a_2a_3-a_2^3-a_0a_3^2-a_1^2a_4,$$
\end{center}
 which is, up to a constant factor, just $\det C(f)_2$ from Example \ref{ex1}, is called classically the Hankel determinant, and the induced function on $S_4$ given by $f \mapsto \det C(f)_2$ is, up to a constant factor, the catalecticant from classical Invariant Theory, see \cite{Do}, p. 10. In particular, the catalecticant is invariant with respect to the group $SL_2(\C)$.
Another invariant of the binary form $f$ is given by
$$S(f)=a_0a_4-4a_1a_3+3a_2^2,$$
see \cite{A,Do}. Using these two invariants, one defines
$$j(f)= \frac{S(f)^3}{ S(f)^3-27T(f)^2}.$$
It is known that two binary quartics $f$ and $f'$, without multiple factors and regarded as points in
$\PP(S_4)$, are in the same $SL_2(\C)$-orbits if and only if
\begin{equation}
\label{eq4}
j(f)=j(f').
\end{equation}

We have the following.
\begin{thm}
\label{thm1}
The Waring rank of a quartic binary form $f$ having no multiple factor is
$2$ if and only if $j(f)=1$.
Otherwise $\rank f =3.$
\end{thm}
\proof
Up to projective equivalence a quartic binary form $f$ having no multiple factor can be written as
\begin{equation}
\label{eq41}
f=xy(x+y)(x+ty),
\end{equation}
with $t \in \C \setminus \{0,1\}$.
Note that one has, using the above formulas,
$$j(f)=\frac{4}{27}\frac{(t^2-t+1)^3}{t^2(t-1)^2}.$$
Since $f$ having no multiple factor, it is clear that $\rank f \geq 2$ by Lemma \ref{lemS}. On the other hand,
Theorem \ref{thmAH} implies that 
$$\rank f \leq d-1=3.$$
Moreover, Lemma \ref{lemS} (3) and our hypothesis that $t \in \C \setminus \{0,1\}$, implies that $\rank f = 2$ if and only if $Ann(f)_2 \ne 0$. Using the formula for the catalecticant $C(f)_2$ given in Example \ref{ex1} and the formula for $f$ in \eqref{eq41}, it follows that
$\det C(f)_2=0$ exactly for $t\in \{-1, {1 \over 2}, 2\}.$ For all these three values of $t$ we get $j(f)=1$.

\endproof

\begin{cor}
\label{cor1}
The quartic binary forms $f$ having no multiple factor and with Waring rank $2$ form a single $SL_2(\C)$-orbit in $\PP(S_4)$. More precisely, $\rank f=2$ if and only if the four roots of $f$, regarded as points in the  projective line $\PP^1$, have a harmonic cross-ratio. 
\end{cor}
\proof
It is known that $j(f)=1$ corresponds exactly to the case when
the four roots of $f$, regarded as points in the  projective line $\PP^1$, have a harmonic cross-ratio. Recall also our remark related to \eqref{eq10} in the Introduction.
\endproof

\begin{rk}
\label{rk2} 
The fact that a binary quartic has Waring rank 2 when the catalecticant $C(f)_2$ vanish and the relation to harmonic cross-ration is stated as a remark in \cite{Ol}, see middle of page 29,
with a reference to an exercise in Gurevich book \cite{Gu}, namely
Exercise 25.7. We leave the interested reader to compare the two different approaches and to notice the distinct terminology used by various authors.
\end{rk}

\section{Binary quintics and the golden ratio}

In this section we investigate the Waring rank of binary forms of degree $5$ having no multiple factor. Up to projective equivalence such a form $f$ can be written as
\begin{equation}
\label{eq5}
f_{s,t}=xy(x+y)(x+sy)(x+ty)=xy(x+y)(x^2+Sxy+Py^2)=f_{S,P},
\end{equation}
with $s,t \in \C \setminus \{0,1\}$ and $s \ne t$. Here $S=s+t$ and $P=st$.
Recall that the golden ratio
$$ \varphi^+={1+\sqrt 5 \over 2}$$
is the positive root of the equation $z^2-z-1=0$. We have the following.
\begin{thm}
\label{thm2}
The Waring rank of the quintic binary form $f_{s,t}$ having no multiple factor is
$2$ if and only if the pair $(s,t)$ is one of the following 12 pairs
$$(\varphi^{\pm}, 1+\varphi^{\pm}), (1+\varphi^{\pm},\varphi^{\pm}),
(-\varphi^{\pm}, 1+\varphi^{\pm}), (1+\varphi^{\pm}, -\varphi^{\pm}),
  (-1+\varphi^{\pm},\varphi^{\pm}), 
(\varphi^{\pm}, -1+\varphi^{\pm}), $$
where $\varphi^{\pm}$ are the two roots of the equation $z^2-z-1=0$.
Otherwise $3 \leq \rank f \leq 4.$ More precisely, the rank of the form 
$ f_{S,P}$ is 4  exactly when the pair $(S,P)$ is a zero of the  polynomial 
$$ \Delta(S,P)=-4S^{12}+12S^{11}P+S^{10}P^2-22S^9P^3+S^8P^4+12S^7P^5-4S^6P^6+ $$
$$+12S^{11}+30S^{10}P-202S^9P^2+84S^8P^3+292S^7P^4-78S^6P^5-102S^5P^6+36S^4P^7+$$
$$+S^{10}-202S^9P+190S^8P^2+1176S^7P^3-1198S^6P^4-1234S^5P^5+666S^4P^6+188S^3P^7-$$
$$-83S^2P^8-22S^9+84S^8P+1176S^7P^2-2640S^6P^3-2264S^5P^4+5392S^4P^5+1236S^3P^6-$$
$$-1532S^2P^7+130SP^8+8P^9 
+S^8+292S^7P-1198S^6P^2-2264S^5P^3+9312S^4P^4-$$
$$-1924S^3P^5-8100S^2P^6+1860SP^7+77P^8 
+12S^7-78S^6P-1234S^5P^2+5392S^4P^3-$$
$$-1924S^3P^4-10570S^2P^5+8010SP^6+120P^7-4S^6-102S^5P+666S^4P^2+1236S^3P^3-$$
$$-8100S^2P^4+8010SP^5-410P^6+36S^4P+188S^3P^2-1532S^2P^3+$$
$$+1860SP^4+120P^5- 83S^2P^2+130SP^3+77P^4+8P^3.$$

%%\label{eqD}

\end{thm}

\proof
As in the proof above, we see that $\rank f \geq 2$ and the equality holds if and only if the catalecticant $C(f)_2$ has not maximal rank 3.
A direct computation shows that
\begin{center}
$$ C(f)_2=\left(
  \begin{array}{ccccccc}
     0 &  4& 2(s+t+1)  \\
     12 &6(s+t+1) &  6(s+t+st) \\
    6(s+t+1)&  6(s+t+st)& 12st  \\
    2(s+t+st) & 4st & 0 \\
   \end{array}
\right).$$
\end{center}
Using the software SINGULAR, we see that the ideal of 3-minors of this matrix has as zero set exactly the 12 pairs $(s,t)$ listed above. 
Assume now that the catalecticant $C(f)_2$ has  maximal rank 3, which implies that $g_1$, the generator of $Ann(f)$ of minimal degree has degree $d_1=3$. It follows that in this case
$\rank (f) \geq 3$. If $g_1=aX^3+3bX^2Y+3cXY^2+dY^3$ is in $Ann(f)_3$, it follows that $(a,3b,3c,d)$ is in the kernel of the matrix $C(f)_3$, and hence in the kernel of the matrix
\begin{center}
$$ \left(
  \begin{array}{ccccccc}
     0 &  2& 1+S &S+P  \\
     2 & 1+S &  S+P& 2P \\
    1+S&  S+P& 2P& 0  \\
   
   \end{array}
\right),$$
\end{center}
obtained by dividing the rows in $C(f)_3$ by 6,12 and 6. Let $m_i$ be the determinant of the $3 \times 3$ matrix obtained from this matrix by deleting the $i$-th column. Then, since the matrix $C(f)_3$ is essentially the transpose of the matrix $C(f)_2$, we know that at least one of the minors $m_i$ is not zero. It follows that one can take
$$(a,3b,3c,d)=(m_1,-m_2,m_3,-m_4),$$
and in this way $a,b,c,d$ become polynomials in $S,P$. 
We define $\al=ac-b^2$, $\be=ad-bc$, $\gamma=bd-c^2$ and
$$\Delta(S,P)=\be^2-4\al \gamma.$$
Then it is known that the binary cubic form $g_1$ has no multiple factors if and only if $\Delta(S,P) \ne0$. In this case $\rank f_{S,P}=3$, and otherwise
$\rank f_{S,P}=4$. This follows from Sylvester's Theorem \ref{thmS}, recalling that $d_1+d_2=d+2=7$ in our case.
\endproof

\begin{rk}
\label{rk2.1} 
Note that a binary form of Waring rank two has necessarily distinct factors, see
\cite[Corollary 4.1.1]{BM}. 
\end{rk}
To the quintic form $f_{s,t}$ above we can associated 5 binary quartic forms without multiple factors, namely
$h_1=f_{s,t}/x$, $h_2=f_{s,t}/y$, $h_3=f_{s,t}/(x+y)$, $h_4=f_{s,t}/(x+sy)$ and $h_5=f_{s,t}/(x+ty)$. It is known that the $SL_2(\C)$-orbit  of $f_{s,t}$ in $\PP(S_5)$ is determined by the unordered list of 5 complex numbers
$$j(f_{s,t}):=((j(h_1),j(h_2),j(h_3),j(h_4),j(h_5)),$$
see \cite[Theorem 13]{A}.
\begin{cor}
\label{cor2}
The quintic binary forms $f$ with Waring rank
$2$ form a single $SL_2(\C)$-orbit in $\PP(S_5)$. More precisely, $\rank f=2$ if and only if $f$ has distinct factors and
$$j(f)=\left(\frac{2^5}{3^3}, \frac{2^5}{3^3},\frac{2^5}{3^3},\frac{2^5}{3^3}, \frac{2^5}{3^3}\right).$$

\end{cor}
\proof
The fact that the 12 pairs $(s,t)$ listed in Theorem \ref{thm2} give rise to a single $SL_2(\C)$-orbit in $\PP(S_5)$ follows from our general remark related to \eqref{eq10} in Introduction.
 A direct computation shows that
$$j(f_{\varphi^{+}, 1+\varphi^{+}})= \left(\frac{2^5}{3^3}, \frac{2^5}{3^3},\frac{2^5}{3^3},\frac{2^5}{3^3}, \frac{2^5}{3^3}\right),  $$
and this completes the proof.
\endproof

\begin{rk}
\label{rk2.5} 
It is shown in \cite[Theorem 4.14]{BM} that there are exactly 
$$N_d={d-1 \choose 2}$$
 distinct forms in $\PP(S_d)$ which are multiple of a fixed cubic form $c\in S_3$ with distinct factors, say $c=xy(x+y)$. For $d=4$ we get $N_4=3$, which explains why we get 3 values for $t$ in the proof of Theorem \ref{thm1} above. Similarly, for $d=5$ we get $N_5=6$, and the corresponding 6 forms are those listed at the beginning of the proof of Corollary \ref{cor2} above.These specific binary forms are related to the map  $\Gamma$, the dihedral cover for the cubic $xy(x+y)$, see \cite[Definition 4.7]{BM}.
\end{rk}

\begin{ex} 
\label{ex3}
Consider the quintic binary form 
$$f=xy(x+y)(x^2+y^2)$$
corresponding to the case $t=i$, $s=-i$ with $i^2=-1.$
The corresponding pair $(S,P)=(s+t,st)$ is now $(0,1)$ and clearly
$\Delta(0,1)=0$. The corresponding form $g_1=(Y-X)(X+Y)^2$ has a multiple factor, and hence 
$$\rank xy(x+y)(x^2+y^2)=4.$$

\end{ex}

\section{On the Waring rank of some ternary forms}

Let $f\in S_d=\C[x,y]_d$ be a binary form of degree $d$ and rank $\rank f \geq 2$, and consider the ternary form $F=zf \in R_{d+1}$, where $R=\C[x,y,z]$.
Assume, using Theorem \ref{thmS},  that $Ann(f)=(g_1,g_2)$ such that 
$$2 \leq d_1 =\deg g_1 \leq d_2=\deg g_2 \text{ and } d_1+d_2=d+2.$$
Then it is clear that $Ann(F)$ in the ring $T=\C[X,Y,Z]$, where $Z$ corresponds to $\partial_z$, it is given by
$$(g_1,g_2,Z^2).$$
The following result is a special case of  \cite[Theorem 4.14]{Car2+}. We include a proof, essentially the same as the proof given in \cite[Theorem 4.14]{Car2+}, just for the reader's convenience.
\begin{thm}
\label{thm3}
The Waring rank of the ternary form $F=zf(x,y)$ is exactly $d_1d_2$, and all the linear forms $\ell_j$ occurring in a minimal length Waring decomposition \eqref{eq1} have the forms $\ell_j=a_jx+b_by+c_jz$ with $c_j \ne 0$ for all $j=1,...,r=d_1d_2$.
\end{thm}

\proof Since $Ann(f)$ is a complete intersection, it follows that $g_2$ can be chosen without multiple factors. With such a choice, we claim that the ideal 
$$I=(g_1(X,Y)+Z^{d_1},g_2(X,Y)) \subset Ann(F) \subset T$$
is a smooth complete intersection $V$, containing $d_1d_2$ simple points in $\PP^2$. Take a point $(p:q:r) \in \PP^2$ in the zero set of this ideal $I$. Note that $r \ne 0$, since the equations 
$$g_1(X,Y)=g_2(X,Y)=0$$
have only the trivial solution $(p,q)=(0,0)$ in $\C^2$. Hence we can take $r=1$
and compute the Jacobian matrix of the mapping $(g_1(X,Y)+Z^{d_1},g_2(X,Y))$
at the point $(p:q:1)$. This matrix has rank 2, due to the fact that
$g_2$ was supposed without multiple factors, and hence
$g_2(p,q)=0$ implies that the gradient of $g_2$ at $(p,q)$ is non-zero.
It follows that 
$$I(V)=I \subset Ann(F).$$
It is known that the Waring rank $\rank F$ is the minimal cardinality of a finite set $W$ in $\PP^2$ such that $I(W) \subset Ann(F)$. 
The set $V$ constructed above shows that this minimal number is $
|W|\leq d_1d_2$.
Let $W'=W \setminus L$, where $L$ is the line $Z=0$. Then one has the following
\begin{enumerate}
\item The cardinality $|W'|$ is equal to the Hilbert polynomial of the quotient $T/I(W')$, which is a constant;

\item Since $Z$ is not a zero-divisor on $T/I(W')$, the above Hilbert polynomial is just the $\C$-dimension of the Artinian algebra $T/(I(W')+(Z)$;

\item Since $I(W) \subset Ann(F)$, we have
$$I(W')=I(W):(Z) \subset Ann(F) :(Z)=(g_1,g_2,Z).$$
\end{enumerate}
It follows that
$$|W'| =\dim \frac{T}{I(W')+(Z)} \geq \dim \frac{T}{(Ann(F) :(Z))+(Z)}=
\dim \frac{T}{(g_1,g_2,Z)}.$$
On the other hand, we have
$$\dim \frac{T}{(g_1,g_2,Z)}=\dim  \frac{Q}{(g_1,g_2)}=d_1d_2,$$
and this proves our claim.
\endproof

\begin{ex} 
\label{ex4}
In the previous sections, we have given examples of quartic binary forms $f$ (respectively quintic binary forms $f$) without multiple factors and such that $(d_1,d_2)=(2,4)$ and $(d_1,d_2)=(3,3)$ for quartic forms, and respectively $(d_1,d_2)=(2,5)$ and $(d_1,d_2)=(3,4)$ for quintic forms. Note that the associated line arrangements in $\PP^2$, namely
$$\A(F): zf(x,y)=0,$$
have a very simple combinatorics, namely a pencil of 4 or 5 lines through a common point, plus a transversal line.
However,  the Waring rank of $F$ can be
8 or 9 for a quartic form $f$, and 10 or 12 for a quintic form $f$.
In particular, the combinatorics of the line arrangement $\A(F)$ cannot determine the Waring rank of the defining equation $F=zf$ of the line arrangement.
\end{ex} 

\begin{rk}
\label{rk3} 
Note the following analog of Lemma 2.2 (3) above. A line arrangement
$\A: F(x,y,z)=0$ in $\PP^2$ satisfies $Ann(F)_1 =0$ and $Ann(F)=\C\ell^2$ for some linear form $\ell \in R_1$ if and only if $\A$ has the same combinatorics as the line arrangements considered in Example \ref{ex4}, namely a pencil of $d$ lines through a point of $\PP^2$ plus a transversal line.
Note that any such arrangement, regarded as a central arrangement in $\C^3$ is free with exponents $(e_1,e_2,e_3)=(1,1,d-1)$, which are the degrees of a basis for the free $R$-module of derivations $D(\A)$,
see for instance 
\cite[Chapter 8]{DHA} for generalities on free arrangements.
On the other hand, the generators of the ideal $Ann(F) \subset T$, as we have seen above, have degrees $(2,d_1,d_2)$ with $d_1+d_2=d+2$.
It does not seem to be a simple relation between the module of derivations $D(\A)$ and the ideal $Ann(F)$, even in this simple situation.
\end{rk}


\begin{thebibliography}{00}

%\bibitem{Abd}  N. Abdallah,  On Hodge theory of singular plane curves. Canad. Math. Bull. 59 (2016),  449--460.

\bibitem{A} A. Abdesselam, A computational solution to a question by Beauville
on the invariants of the binary quintic, J. Algebra 303 (2006), 771--788.

\bibitem{AH}  J. Alexander, A. Hirschowitz, Polynomial interpolation in several variables, J. Algebraic Geom. 4 (1995),  201--222.

%\bibitem{Abe18}  T. Abe, Plus-one generated and next to free arrangements of hyperplanes, arXiv:1808.04697.

%\bibitem{AD}  T. Abe, A. Dimca, On the splitting types of bundles of logarithmic vector fields along plane curves, Internat. J. Math. 29 (2018), no. 8, 1850055, 20 pp.

%\bibitem{B+} E. Artal Bartolo, L. Gorrochategui, I. Luengo, A. Melle-Hern\' andez, On some conjectures about free and nearly free divisors, in: {\it Singularities and Computer Algebra, Festschrift for Gert-Martin Greuel on the Occasion of his 70th Birthday}, pp. 1--19, Springer (2017)



\bibitem{Guide} A. Bernardi, E. Carlini, M. V. Catalisano, A. Gimigliano, A. Oneto, The Hitchhiker guide to: Secant Varieties and Tensor Decomposition, arXiv:1812.10267.

\bibitem{BM} L. Brustenga i Moncus\' i, S. K. Masuti, On the Waring rank of binary forms: the binomial formula and a dihedral cover of rank two forms, arXiv:1901.08320.

\bibitem{BT} J. Buczy\' nski, Z. Teitler, Some examples of forms of high rank, Collectanea Mathematica 67(2016), 431--441.



\bibitem{Car1+} E. Carlini, M.V. Catalisano, A. Oneto,  Waring loci and the Strassen conjecture. Adv. Math. 314(2017), 630--662.


\bibitem{Car2+} E. Carlini, M. V. Catalisano, L. Chiantini, A. V. Geramita, Y. Woo, 
Symmetric tensors: rank, Strassen’s conjecture and  e-computability,
Ann.  Scuola Normale Sup. Pisa. 18 (2018), 363--390.

%\bibitem{Cay} A. Cayley,  A Memoir on Cubic Surfaces. Philos. Trans. R. Soc. Lond., Ser. A 159 (1869), 231--326.

%\bibitem{CD} A.~D.~R.~Choudary, A.~Dimca,  Koszul complexes and hypersurface singularities,  Proc. Amer. Math. Soc. 121(1994), 1009--1016. 

\bibitem{Ci} C. Ciliberto, Geometric aspects of polynomial interpolation in more variables and of Waring’s problem. In: European Congress of Mathematics, Barcelona 2000, pages 289--316. Springer, 2001.

\bibitem{Co} CoCoA-5 (15 Sept 2014): a system for doing Computations in Commutative Algebra,
available at http://cocoa.dima.unige.it

\bibitem{CS} G. Comas and M. Seiguer, On the rank of a binary form. Foundations of Computational
Mathematics, 11(1)52011), 65--78.

\bibitem{CGLM} P. Comon, G. Golub, L.H. Lim, B. Mourrain, Symmetric tensors and symmetric tensor
rank,  SIAM Journal on Matrix Analysis and Applications, 30(2008),1254--1279.

%\bibitem{CHMN} D. Cook II, B. Harbourne, J. Migliore, U. Nagel, Line arrangements and configuration of points with an unusual geometric property, arXiv: 1602.02300.

\bibitem
{Sing} { W. Decker, G.-M. Greuel, G. Pfister \and H. Sch{\"o}nemann.} \newblock {\sc Singular} {4-0-1} --- {A} computer algebra system for polynomial computations, available at {http://www.singular.uni-kl.de} (2014).


%\bibitem{DBull}  A. Dimca, Syzygies of Jacobian ideals and defects of linear systems, Bull. Math. Soc. Sci. Math. Roumanie Tome 56(104) No. 2, 2013, 191--203.


%\bibitem{DRCS}  A. Dimca,   {\em Topics on Real and Complex Singularities}, Vieweg Advanced Lecture in Mathematics, Friedr. Vieweg und Sohn, Braunschweig, 1987.

%\bibitem{DSTH}  A. Dimca,   {\em Singularities and Topology of Hypersurfaces} , Universitext, Springer Verlag, New York, 1992.

%\bibitem{DST}  A. Dimca,   {\em Sheaves in Topology} , Universitext, Springer, 2004.

\bibitem{DHA}  A. Dimca,   {\em Hyperplane Arrangements: An Introduction}, Universitext, Springer, 2017


%\bibitem{Dmax}  A. Dimca, Freeness versus maximal global Tjurina number for plane curves,  Math. Proc. Cambridge Phil. Soc.  163 (2017), 161--172.

%\bibitem{Drcc} A. Dimca, On rational cuspidal plane curves, and the local cohomology of Jacobian rings, arXiv:1707.05258.  to appear in Commentarii Mathematici Helvetici.


%\bibitem{DIM} A. Dimca, D. Ibadula, A. M\u acinic, Numerical invariants and moduli spaces for line arrangements, arXiv:1609.06551.


%\bibitem{DPop} A. Dimca, D. Popescu, Hilbert series and Lefschetz properties of dimension one almost complete intersections, Comm. Algebra 44 (2016), 4467--4482.




%\bibitem{DS14} A. Dimca, E. Sernesi,  Syzygies and logarithmic vector fields along plane curves, Journal de l'\'Ecole polytechnique-Math\'ematiques 1(2014), 247-267.

%\bibitem{Camb2012} A. Dimca, G. Sticlaru, Chebyshev curves, free resolutions and rational curve arrangements, Math. Proc. Cambridge Phil. Soc. 153  (2012), 385­--397.


%\bibitem{DStExpo} A. Dimca, G. Sticlaru, On the exponents of free and nearly free projective plane curves, Rev. Mat. Complut. 30(2017), 259--268.

%\bibitem{DStFD} A. Dimca, G. Sticlaru, Free divisors and rational cuspidal plane curves, Math. Res. Lett. 24(2017), 1023--1042.

%\bibitem{DStNF} A. Dimca, G. Sticlaru, Nearly free divisors and rational cuspidal curves, arXiv:1505.00666.


%\bibitem{DStRIMS} A. Dimca, G. Sticlaru, Free and nearly free curves vs. rational cuspidal plane curves, Publ. RIMS Kyoto Univ. 54 (2018), 163--179.

%\bibitem{DStMos} A. Dimca, G. Sticlaru, On the freeness of rational cuspidal plane curves, arXiv:1802.06688, to appear in Moscow Math.J.
  


%\bibitem{DStFS} A. Dimca, G. Sticlaru, Free and nearly free surfaces in $\PP^3$, arXiv:1507.03450v3.

%\bibitem{DStexpo} A. Dimca, G. Sticlaru, On the exponents of free and nearly free projective plane curves, arXiv:1511.08938.

%\bibitem{DStSat} A. Dimca, G. Sticlaru, Saturation of Jacobian ideals: some applications to nearly free curves, line arrangements and rational cuspidal plane curves, arXiv: 1711.02595v2.

%\bibitem{DStJump} A. Dimca, G. Sticlaru, On the jumping lines of bundles of logarithmic vector fields along plane curves, arXiv: 1804.06349.

\bibitem{Do} I. Dolgachev, Lectures on Invariant Theory, London Mathematical Society Lecture Note Series,  Cambridge: Cambridge University Press, 2003. 

%\bibitem{duPCTC} A.A. du Plessis,  C.T.C. Wall, Application of the theory of the discriminant to highly singular plane curves, Math. Proc. Cambridge Phil. Soc.,  126(1999), 259-266. 

%\bibitem{duPCTC2} A.A. du Plessis and C.T.C. Wall, Curves in $P^2(\C)$ with 1-dimensional symmetry, Revista Mat Complutense 12 (1999), 117--132.

%\bibitem{Eis} { D. Eisenbud}, \emph{The Geometry of Syzygies: A Second Course in Algebraic Geometry and Commutative Algebra}, Graduate Texts in Mathematics, Vol. 229, Springer 2005. 

%\bibitem{FV1}  D. Faenzi, J. Vall\`es, {Logarithmic bundles and line arrangements, an approach via the standard construction}, {J. London.Math.Soc.} {90} (2014), {675--694}.


%\bibitem{FZ} H. Flenner,  M. Zaidenberg, On a class of rational plane curves, Manuscripta Math. 89 (1996), 439--459.

%\bibitem{FLMN} J. Fernandez de Bobadilla, I. Luengo, A. Melle-Hernandez, A. Nemethi, Classification of rational unicuspidal projective curves whose singularities have one Puiseux pair,
%in: {\it Real and Complex Singularities Sao Carlos 2004}, Trends in Mathematics, Birkhauser 2006, pp. 31-45.

\bibitem{FOS} R. Fr\" oberg, G. Ottaviani,  B. Shapiro, On the Waring problem for polynomial rings Proceedings
of the National Academy of Sciences, 109 (2012), 5600--5602.

\bibitem{FLOS} R. Fr\" oberg,  S. Lundqvist, 
A. Oneto, 
B. Shapiro, 
Algebraic stories from one and from the other pockets, Arnold Math. J. 4 (2018),  137--160.

%\bibitem{Ge} A.V. Geramita. Inverse systems of fat points: Waring’s problem, secant varieties of Veronese varieties and parameter spaces for Gorenstein ideals. In: The Curves Seminar at Queen’s, volume 10, pages 2--114, 1996.

%\bibitem{GV} B. Guerville-Ball\'e, J. Viu-Sos, Combinatorics of line arrangements and dynamics of polynomial vector fields, arXiv:1412.0137.

\bibitem{Gu} G. B. Gurevich, Foundations of the Theory of Algebraic Invariants, P. Noordhoff,
1964.

%\bibitem{H+1} T. Harima, J. Migliore, U. Nagel and J. Watanabe, { The weak and strong
%Lefschetz properties for artinian K-algebras}, J. Algebra \textbf{262} (2003), 99--126.

%\bibitem{H+2} T. Harima,  T. Maeno, H.  Morita, Y.  Numata, A.  Wachi, J.  Watanabe, {\it The Lefschetz properties}, Lecture Notes in Mathematics 2080, Springer, Heidelberg, 2013.

%\bibitem{HS} S. H. Hassanzadeh, A. Simis, Plane Cremona maps: Saturation and regularity of the base ideal, J. Algebra 371 (2012), 620--652.

%\bibitem{I} A. Iarrobino, Compressed Algebras: Artin Algebras Having Given Socle,  Degrees and Maximal Length, Trans. Amer. Math. Soc., 285 (1984), 337--378.

\bibitem{IK} A. Iarrobino, V. Kanev, {\it Power Sums, Gorenstein Algebras, and Determinantal Loci}, Springer Lecture Notes 1721, 1999.

%\bibitem{IT} N. Ilten, Z.Teitler,  Product ranks of the $3 \times 3$ determinant and permanent. Canad. Math. Bull. 59 (2016), 311--319.


%\bibitem{IS} N. Ilten, H. S\"u\ss, Fano schemes for generic sums of products of linear forms, arXiv:1610.06770.



%\bibitem{Kl} S.L. Kleiman,  The enumerative theory of singularities, in: Real and Complex Singularities (Oslo 1976), Sijthoff and Noordhoff, Amsterdam 1977, pp. 297-- 396.

\bibitem{L1} J.M. Landsberg, {\it Tensors: Geometry and Applications}, Graduate Studies in Mathematics vol.128, American Mathematical Soc., 2012.

\bibitem{LT} J.M. Landsberg, Z. Teitler, On the ranks and border ranks of symmetric tensors, Found. Comput. Math. 10(3) (2010) 339--366.



%\bibitem{IG} G. Ilardi, { Jacobian ideals, Arrangement and Lefschetz properties}, J. Algebra  \textbf{508} (2018), 418--430. 

%\bibitem{JF} G. Jiang, J. Feng, An algorithm to produce the structure sequence of an arrangement, Singularities and complex geometry, Q. Lu, S.S.T Yau, A. Libgober, eds., AMS/IP Studies in Advanced Mathematics, 5 (1997), 83--92.

%\bibitem{Lin} N. Lindner, Cuspidal plane curves of degree 12 and their Alexander polynomials,
 %(64 pages, Master thesis, Humboldt Universit\"at Berlin, 2012)
 
% \bibitem{MaVa} S. Marchesi, J. Vall\` es, Nearly free curves and arrangements: a vector bundle point of view, arXiv:1712.04867.


%\bibitem{Moe} T. K.  Moe, Rational Cuspidal Curves, arXiv:1511.02691 (139 pages, Master thesis 2008)

\bibitem{MO} B. Mourrain, A. Oneto,  On minimal decompositions of low rank symmetric tensors, arXiv:1805.11940.

\bibitem{Ol} P. Olver, Classical Invariant Theory, London Mathematical Society Student
Texts 44, Cambridge Univ. Press, 1999.

\bibitem{O} A. Oneto, Waring type problems for polynomials, Doctoral Thesis in Mathematics at Stockholm University, Sweden, 2016.


%\bibitem{OT} P. Orlik and H. Terao, {\em Arrangements of Hyperplanes,} Springer-Verlag, Berlin Heidelberg New York, 1992.

%\bibitem{Pion} J. Piontkowski, On the number of cusps of rational cuspidal plane
%curves, Experiment. Math., 16(2007), 251-255.

%\bibitem{KS} K. Saito, Quasihomogene isolierte Singularit\"aten von Hyperfl\"achen, Invent. Math., 14 (1971), 123--142.

%\bibitem{KS} K. Saito, Theory of logarithmic differential forms and logarithmic vector fields, J. Fac. Sci. Univ. Tokyo Sect. IA Math. 27 (1980), no. 2, 265-291.

%\bibitem{SaTo} F. Sakai, K. Tono, Rational cuspidal curves of type (d,d-2) with one or two cusps, Osaka J. Math. 37(2000), 405-415.

%\bibitem{Sch0} H. Schenck, Elementary modifications and line configurations in $P^2$,  Comm. Math. Helv. 78 (2003), 447--462.

%\bibitem{SchToh} H. Schenck, S. O. Toh\u aneanu, 
%Freeness of conic-line arrangements in $P^2$,  Comm. Math. Helv. 84 (2009), 235--258.

%\bibitem{Schenck} H. Schenck, Hyperplane arrangements: computations and conjectures. Arrangements of hyperplanes--Sapporo 2009, 323--358, Adv. Stud. Pure Math., 62, Math. Soc. Japan, Tokyo, 2012.


%\bibitem{Se} E. Sernesi,  The local cohomology of the jacobian ring, {Documenta Mathematica},  19 (2014), 541-565. 

%\bibitem{Sim1} A. Simis,  Differential idealizers and algebraic free divisors, in:{\it Commutative Algebra: Geometric, Homological,Combinatorial and  Computational Aspects}, Lecture Notes in Pure and Applied Mathematics (Eds. A. Corso, P. Gimenez, M. V. Pinto and S. Zarzuela), Chapman \& Hall, Volume 244 (2005),211--226.


%\bibitem{Sim2} A. Simis, The depth of the Jacobian ring of a homogeneous polynomial in three variables, Proc. Amer. Math. Soc., 134 (2006), 1591--1598.

%\bibitem{ST} A. Simis, S.O. Toh\u aneanu, Homology of homogeneous divisors, Israel J. Math. 200 (2014), 449-487.

%\bibitem{SW} D. van Straten, T. Warmt,  Gorenstein duality for one-dimensional almost complete intersections--with an application to non-isolated real singularities, Math. Proc.Cambridge Phil. Soc.158(2015), 249--268.

\bibitem{Sy} J.J. Sylvester. Lx. on a remarkable discovery in the theory of canonical forms and of hyperdeterminants.
The London, Edinburgh, and Dublin Philosophical Magazine and Journal of
Science, 2(12)(1851), 391--410.

%\bibitem{Wa} U. Walther, The Jacobian module, the Milnor fiber, and the $D$-module generated by $f^s$,  Invent. Math. 207 (2017),1239--1287.

%\bibitem{ZW} Z. Wang, On homogeneous polynomials determined by their partial derivatives, to appear in
Canad. Math. Bull. 2019. DOI: https://doi.org/10.4153/S0008439519000419.

%\bibitem{Yo} M. Yoshinaga, Freeness of hyperplane arrangements and related topics, Annales de la Facult\'e des Sciences de Toulouse, vol. 23 no. 2 (2014), 483-512.

%\bibitem{Zi} G. Ziegler, Combinatorial construction of logarithmic differential forms, Adv. Math. 76 (1989), 116-154.

\end{thebibliography}
\end{document}